\documentclass[12pt,letterpaper]{article}

\usepackage{graphicx}
\usepackage{psfrag}
\usepackage{url}
\usepackage{epsfig}
\usepackage{fancybox}
\usepackage{geometry}
\usepackage[nosort, nomove]{cite}
\usepackage{amsmath,amssymb,mathrsfs,amsthm}
\usepackage{subfigure}
\usepackage{rotating}
\usepackage{color}
\usepackage{calc}
\usepackage{setspace}

\newcommand{\ld}{\lambda}

\newcommand{\rset}{\mathbb{R}}

\newcommand{\xb}{\textbf{x}}

\newcommand{\bmtrx}{\left[\begin{array}}
\newcommand{\emtrx}{\end{array}\right]}

\newtheorem{theorem}{Theorem}[section]
\newtheorem{lemma}{Lemma}[section]
\newtheorem{proposition}{Proposition}[section]
\newtheorem{corollary}{Corollary}[section]

\theoremstyle{definition}
\newtheorem{remark}{Remark}[section]
\newtheorem{definition}{Definition}[section]
\newtheorem{algorithm}{Algorithm}[section]
\newtheorem{assumption}{Assumption}[section]

\title{An Interior-Point Lagrangian Decomposition Method for Separable Convex Optimization}
\author{I. Necoara$^{1,3}$ and  J.A.K. Suykens$^{2,3}$}%
\date{September 2008}
\date{Communicated by   D. Q. Mayne}

\graphicspath{{./}{figures/}{pictures/}}

\begin{document}
\allowdisplaybreaks

\maketitle

\footnotetext[1]{Politehnica University of Bucharest,  Automatic
Control and Systems Engineering  Department, 060042 Bucharest,
Romania  \& Katholieke Universiteit Leuven, Department of Electrical
Engineering (ESAT),  B--3001 Leuven--Heverlee, Belgium.}

\footnotetext[2]{Katholieke Universiteit Leuven, Department of
Electrical Engineering (ESAT), B--3001 Leuven--Heverlee, Belgium.}

\footnotetext[3]{We  acknowledge  financial support from  Flemish
Government:  FWO projects G.0226.06, G.0302.07.}


\setcounter{footnote}{4}
\newpage

\noindent {\bf Abstract.} In this paper, we propose a distributed
algorithm for solving large-scale separable convex problems using
Lagrangian dual decomposition and the interior-point framework. By
adding self-concordant barrier terms to the ordinary Lagrangian, we
prove under mild assumptions  that the corresponding family of
augmented dual functions is self-concordant. This makes it possible
to efficiently use the Newton method for tracing the central path.
We show that the new algorithm is globally convergent and  highly
parallelizable and thus it is suitable for solving large-scale
separable convex problems.


\newpage

\section{Introduction}

Can self-concordance and interior-point methods be incorporated into
a Lagrangian dual decomposition  framework? This paper presents a
decomposition algorithm   that incorporates the interior-point
method into augmented Lagrangian decomposition technique for solving
large-scale separable convex problems.  Separable convex
 problems, i.e. optimization problems with a separable convex objective function
  but with coupling
constraints,  arise in many fields:  networks (communication
networks, multicommodity network flows) \cite{XiaBoy:04,GonSar:03},
process system engineering (e.g. distributed model predictive
control) \cite{VenRaw:07,Nec:08}, stochastic programming
\cite{Zha:05}, etc. There has been considerable interest in parallel
and distributed computation methods for solving this type of
structured optimization problems and many methods have been
proposed: dual subgradient methods \cite{BerTsi:89,Lem:06},
alternating direction methods \cite{BerTsi:89,KonLeo:96}, proximal
method of multipliers \cite{CheTeb:94},  proximal center method
\cite{Nec:08}, interior-point  based methods
\cite{GonSar:03,Zha:05,KojMeg:93,KorPot:91,Tse:92,Zhu:92,HegOsb:01},
etc.

The methods mentioned above belong to the class of augmented
Lagrangian or multiplier methods \cite{BerTsi:89}, i.e. they can be
viewed as techniques for maximizing an augmented dual function. For
example in the  alternating direction   method a quadratic penalty
term is added to the standard Lagrangian to obtain a smooth dual
function and then using a steepest ascent update for the
multipliers. However, the quadratic term destroys the separability
of the given problem. 
Moreover,  the performance of these  methods is very sensitive to
the variations of their parameters and some rules for choosing these
parameters were given e.g. in \cite{KonLeo:96,MieMos:72}. In the
proximal center method \cite{Nec:08} we use smoothing techniques in
order to obtain a well-behaved Lagrangian, i.e.  we add a separable
strongly convex term to the ordinary Lagrangian. This  technique
leads to a smooth dual function, i.e. with Lipschitz continuous
gradient, which preserves separability  of the problem, the
corresponding parameter is selected optimally and moreover the
multipliers are updated using an optimal gradient based scheme. In
\cite{GonSar:03,Zha:05,KojMeg:93,KorPot:91,Tse:92,Zhu:92,HegOsb:01}
interior-point methods are proposed for solving special classes of
separable convex problems with a particular structure of the
coupling/local constraints. In those papers the Newton direction is
used to update the primal variables and/or multipliers obtaining
polynomial-time complexity for the proposed algorithms. \\
In the
present paper we use a similar smoothing technique as in
\cite{Nec:08} in order to obtain a well-behaved augmented dual
function.  Although we  relax the coupling constraints using the
Lagrangian dual  framework as in \cite{Nec:08}, the main difference
here  is that  the smoothing term is a self-concordant barrier,
while in \cite{Nec:08} the main property of the smoothing term was
strong convexity. Therefore, using the properties of self-concordant
functions we show  that the augmented dual function becomes under
mild assumptions also self-concordant. Hence the Newton direction
can be used instead of gradient based directions as it is done in
 most of the augmented Lagrangian methods.
Furthermore, we develop a specialized interior-point method to
maximize the augmented dual function which takes into account the
special structure of our problem. We present a parallel algorithm
for computing the Newton direction of the dual function and we also
prove global convergence of the proposed method.

The main contributions of the paper are the following:\\
(i) We consider a more general model  for separable convex problems
that includes  local equality and inequality constraints, and linear
coupling constraints which generalizes the  models in
\cite{GonSar:03,Zha:05,KojMeg:93,HegOsb:01}.\\
(ii) We derive sufficient conditions for self-concordance of
augmented Lagrangian and
we prove self-concordance for the corresponding family of augmented dual functions.\\
(iii) We provide an interior-point based algorithm for solving the
dual problem with proofs of global convergence and polynomial-time
complexity. \\
(iv) We propose a practical implementation of the algorithm based on
solving approximately the subproblems and on parallel computations
 of the Newton directions.

Note that  item (ii) generalizes the results  of
\cite{Zha:05,KojMeg:93}. 
However, the consideration of general convex problems with local
equality constraints requires new proofs with more involved
arguments in order to prove self-concordance for the family  of dual
functions.

 This paper is organized as follows. In Section
\ref{spf} we formulate the separable convex problem followed by a
brief description of some of the existing decomposition methods for
this problem. The main results  are given in Sections \ref{s_sc} and
\ref{s_pflda}. In Section \ref{s_sc} we show that the augmented
Lagrangian obtained by adding self-concordant barrier terms to the
ordinary Lagrangian forms a self-concordant family of dual
functions. Then an interior-point Lagrangian decomposition algorithm
with polynomial complexity is proposed in Section \ref{s_pflda}. The
new algorithm makes use of the special structure of our problem so
that it is highly parallelizable and it can be  effectively
implemented on parallel processors.   We conclude the paper with
some possible applications.

Throughout the paper  we use the following notations. For a function
$\psi$ with two arguments, scalar parameter $t$ and  decision
variable $x$, i.e. $\psi(t,x)$, we use `` {\small{$\prime$}}'' to
denote  the partial derivative of $\psi(t,x)$ with respect to $t$
and ``$\nabla$'' with respect to $x$: e.g. $\nabla \psi^\prime(t,x)
= \frac{\partial^2}{\partial t  \partial x} \psi(t,x)$. For a
function $\phi$, three times differentiable, i.e. $\phi \in {\cal
C}^3(\text{dom} \!\ \phi)$, $\nabla^3 \phi(x) [h_1,h_2,h_3]$ denotes
the third differential of $\phi$ at $x$ along   directions $h_1,
h_2$ and $h_3$. We use the notation $A \preceq B$ if $B-A$ is
positive semidefinite. We use $D_A$ to denote the block diagonal
matrix having on the main diagonal the matrices $A_1,\cdots,A_N$. We
use $\text{int}(X)$ to denote the interior of a set $X$.

\section{Problem Formulation}
\label{spf}

We consider the following general \textit{separable convex}
optimization problem:
\begin{align}
\label{scp_f}
f^* = & \min_{x_1 \in X_1 \cdots x_N \in X_N}   \sum_{i=1}^N f_i(x_i) \\
\label{scp_B} &\text{s.t.} \; \sum_{i=1}^N B_i x_i = b, \;\; A_i x_i
= a_i  \; \forall i=1 \cdots N,
\end{align}
where $f_i: \rset^{n_i} \to \rset$ are  convex functions, $X_i$ are
closed convex sets, $A_i \in \rset^{m_i \times n_i}$, $B_i \in
\rset^{m \times n_i}$, $a_i \in \rset^{m_i}$ and $b \in \rset^m$.
For simplicity of the exposition we define the vector $x:= [x_1^T
\cdots x_N^T]^T$, the function $f(x):~=~\sum_{i=1}^N f_i(x_i)$, the
set $X:= \prod_{i=1}^N X_i$, the matrix  $B:= [B_1 \cdots B_N]$ and
$n:=\sum_{i=1}^N n_i$.
\begin{remark}
(i) Note that we do not assume strict/strong convexity of any
function~$f_i$.\\
(ii) Coupling \textit{inequality} constraints $\sum_{i=1}^N B_i x_i
\leq b$ can be  included in this framework  by  adding a slack
 variable $x_{N+1}$: $\sum_{i=1}^N B_i x_i + x_{N+1} =
 b$, i.e. $B_{N+1} \!=\! I$ and $X_{N+1} = \rset^m_+$.
\end{remark}
Let   $\langle \cdot,\cdot \rangle$ denote the Euclidian inner
product on  $\rset^n$. By forming the Lagrangian corresponding only
to the coupling linear constraints  (with the  multipliers $\ld \in
\rset^m$), i.e.
 \[ L_0(x,\ld) = f(x) + \langle \ld,B x - b \rangle, \]
we can define the standard dual function
\[ d_0(\ld) = \min_x \{L_0(x,\ld): x_i \in X_i,\; A_i x_i = a_i \; \forall i=1 \cdots N \}.  \]
Since $L_0$ preserves the separability of our problem we can  use
the dual decomposition method \cite{BerTsi:89} by solving in
parallel $N$ minimization problems and then updating the multipliers
in some fashion.  Note that the dual function $d_0$ is concave but,
in general  $d_0$ is not differentiable (e.g.  when $f$ is not
strictly convex). Therefore, for maximizing $d_0$ one has to use
 involved nonsmooth optimization techniques  \cite{BerTsi:89,Lem:06}. From  duality theory
 one knows that if $(x^*, \ld^*)$ is
a saddle point for the  Lagrangian $L_0$, then under appropriate
conditions (constraint qualification), $x^*$ is an optimal solution
for the primal \eqref{scp_f}--\eqref{scp_B} and $\ld^*$ is an
associated dual optimal multiplier for the dual problem: $\max_{\ld
\in \rset^m}
  d_0(\ld)$.

 In order to obtain a smooth dual function we need to use smoothing
 techniques applied to the ordinary Lagrangian $L_0$. One approach
 is the \textit{augmented  Lagrangian} obtained e.g. by adding a quadratic
   penalty term to the  Lagrangian $L_0$: $t \|B x - a\|^2 $. In
the alternating direction method \cite{BerTsi:89,KonLeo:96} the
minimization of the augmented Lagrangian is performed  by
alternating minimization in a Gauss-Seidel fashion followed by a
steepest ascent update  for the multipliers.

In \cite{Nec:08} we proposed
 the \textit{proximal center method} in which  we added to the  standard Lagrangian a smoothing
 term $t \sum_{i=1}^N g_{X_i} (x_i)$, where each function $g_{X_i}$ is
strongly convex  and depends on the set $X_i$ so that the augmented
Lagrangian takes the following form:

\[ L_t^{\text{prox}} (x,\ld) = \sum_{i=1}^N [f_i(x_i) +
t  g_{X_i} (x_i)] + \langle \ld, B x - b \rangle.  \]
Therefore, the
augmented Lagrangian $L_t^{\text{prox}}$ is strongly convex,
preserves separability of the problem like $L_0$ and the associated
augmented dual function
\[ d_t^{\text{prox}} (\ld) = \min_{x} \{ L_t^{\text{prox}}(x,\ld) : x_i \in X_i,
 \; A_i x_i = a_i  \; \forall i=1 \cdots N \} \]
is differentiable and has also a Lipschitz
 continuous gradient.
In \cite{Nec:08} an accelerated gradient based method  is used to
maximize the augmented dual function $d_t^{\text{prox}}$,  while the
corresponding minimization problems are solved in parallel.
Moreover, the smoothing parameter $t$ is selected optimally.

Note that  the methods discussed above use only the gradient
directions of the augmented dual function in order to update the
multipliers. Therefore, in the absence of more conservative
assumptions like strong convexity,  the global convergence rate of
these methods  is slow, in general sub-linear. In this paper we
propose to smoothen the Lagrangian by  adding instead of strongly
convex terms $g_{X_i}$, self-concordant barrier terms $\phi_{X_i}$
associated  with the sets $X_i$, in order to obtain the
self-concordant Lagrangian:
\begin{align}
\label{lip} L_t^{\text{sc}} (x,\ld) = \sum_{i=1}^N [f_i(x_i) + t
\phi_{X_i} (x_i)] + \langle \ld, B x - b \rangle.
\end{align}
In the next section we show, using the theory of self-concordant
barrier functions \cite{NesNem:94,Ren:01}, that for a relatively
large class of convex functions $f_i$ (see also Section
\ref{s_apl}), we can obtain a self-concordant augmented dual
function:

\begin{align}
\label{dip} d^{\text{sc}}(t,\ld) = \min_{x} \{
L_t^{\text{sc}}(x,\ld) : x_i \in \text{int}(X_i),  \; A_i x_i = a_i
\; \forall i=1 \cdots N \}.
\end{align}

This opens the possibility of deriving an
interior-point  method using  Newton  directions  for updating the
multipliers to speed up the convergence rate of the proposed
algorithm.


\section{Sufficient Conditions for Self-Concordance of the Augmented Dual Function}
\label{s_sc}

In this section we derive  sufficient conditions under which the
family of augmented dual functions is self-concordant. A key
property that allows to prove polynomial convergence for barrier
type methods is the property of self-concordance (see Definition
2.1.1 in \cite{NesNem:94}):
\begin{definition}
 A closed convex
 function $\phi$ with open convex domain $X_\phi \subseteq \rset^n$ is called
   $M_\phi$-\textit{self-concordant},
  where $M_\phi \geq 0$, if
$\phi$ is three times  continuously differentiable on $X_\phi$ and
if for all $x \in X_\phi$ and $h \in \rset^n$ we have

\begin{align}
\label{d_sc} \nabla^3 \phi(x) [h,h,h]  \leq  M_\phi \big ( h^T
\nabla^2 \phi(x) h \big )^{3/2}.
\end{align}

A function $\phi$  is called
   $N_\phi$-\textit{self-concordant barrier} for its domain $X_\phi$
   if $\phi$ is 2-self-concordant function and for all
   $x \in X_\phi$ and $h \in \rset^n$ we have
\begin{align}
\label{d_scb} \langle \nabla \phi(x), h \rangle^2 \leq  N_\phi \ h^T
\nabla^2 \phi(x) h.
\end{align}
\end{definition}

 Note that \eqref{d_sc} is equivalent to
(see \cite{NesNem:94}, pp. 14):
\begin{align}
\label{d_sce} | \nabla^3 \phi(x) [h_1,h_2,h_3] | \leq  M_\phi
\prod_{i=1}^3 \big ( h_i^T \nabla^2 \phi(x) h_i \big )^{1/2}.
\end{align}

\noindent Moreover, if  Hessian $\nabla^2 \phi (x)$ is positive
definite, then the inequality \eqref{d_scb} is equivalent to
\begin{align}
\label{d_scb1} \nabla \phi(x)^T \nabla^2 \phi(x)^{-1}  \nabla
\phi(x) \leq N_\phi.
\end{align}

Next lemma provides some basic properties
of self-concordant functions:
\begin{proposition} (\cite{NesNem:94}, pp. 15)
\label{l_scp} Let $\phi$ be an  $M_\phi$-self-concordant function
such that its domain $X_\phi$ does not contain straight lines (i.e.
sets of the form  $\{ x + \alpha u:  \alpha \in \rset \}$, where $x
\in X_\phi$  and $u \not = 0$). Then, the Hessian $\nabla^2 \phi
(x)$ is positive definite for all $x  \in X_\phi$ and $\phi$ is a
barrier function for $X_\phi$. \qed
\end{proposition}
Note that a self-concordant function  which is also a barrier for
its domain is called \textit{strongly self-concordant}.  The next
lemma gives some helpful composition rules for self-concordant
functions.
\begin{lemma}
\label{l_scsc} (i) \cite{NesNem:94} Any linear or convex quadratic
function is
0-self-concordant.\\
(ii) \cite{NesNem:94} Let $\phi_i$ be $M_i$-self concordant and let
$p_i>0$, $i=1,2$.
 Then the function $p_1 \phi_1 + p_2 \phi_2$ is also $M$-self
concordant, where $M=\max \{ M_1/\sqrt{p_1}, M_2 / \sqrt{p_2}\}$.\\
(iii) Let $X_\text{box} = \prod_{i=1}^n [l_i, \ u_i]$ such that $l_i
< u_i$ and $\psi \in {\cal C}^3 (\text{int}(X_\text{box}))$ be
convex. If there exists  $\beta>0$ such that for all $x \in
\text{int}(X_\text{box})$ and $h \in \rset^n$ the following
inequality holds %
\begin{align}
\label{l_box} | \nabla^3 \psi(x) [h,h,h]  | \leq \beta \ h^T
\nabla^2 \psi(x) h \ \sqrt{\sum_{i=1}^n  h_i^2/(u_i-x_i)^2
+h_i^2/(x_i-l_i)^2},
\end{align}
 then $\bar \psi_t(x) = \psi(x) -t \sum_{i=1}^n  \log (u_i-x_i)(x_i - l_i)$ is
$2(1+\beta)/\sqrt{t}$-self concordant.
\end{lemma}
\proof (i) and (ii) can be found in \cite{NesNem:94}, pp. 13.\\
(iii) Denote  $\phi_\text{box}(x) = - \sum_{i=1}^n  \log
(u_i-x_i)(x_i - l_i)$. Note that
\begin{align*}
h^T \nabla^2 \phi_\text{box}(x) h & = \sum_{i=1}^n
  h_i^2/(u_i-x_i)^2 + h_i^2/(x_i-l_i)^2  \\
\nabla^3 \phi_\text{box}(x) [h,h,h] & = 2 \sum_{i=1}^n
h_i^3/(u_i-x_i)^3 - h_i^3/(x_i-l_i)^3
\end{align*}

 and using Cauchy-Schwarz inequality it
follows that $\phi_\text{box}$ is   $2$-self-concordant function on
$\text{int}(X_\text{box})$. Let us denote \[ c= \sqrt{h^T \nabla^2
\psi(x) h} \quad  \text{and} \quad  d= \sqrt{\sum_{i=1}^n
  h_i^2/(u_i-x_i)^2 + h_i^2/(x_i-l_i)^2}. \]  Using  hypothesis \eqref{l_box}  and
  2-self-concordance of $\phi_\text{box}$ we
  have the following   inequalities:
\begin{align*}
| \nabla^3 \bar \psi_t(x) [h,h,h]  | \leq  |\nabla^3 \psi(x) [h,h,h]
| + t |\nabla^3 \phi_\text{box}(x) [h,h,h]| \leq \beta c^2 d + 2t
d^3.
\end{align*}

With some computations we can observe that
\[ (\beta c^2 d + 2t d^3)^2 \leq \frac{4(1+\beta)^2}{t} (c^2 + t
d^2)^3\]  and since \[ h^T \nabla^2 \bar \psi_t(x) h = c^2 + t d^2,
\] the proof is complete. \qed

\noindent  Note that condition \eqref{l_box} is similar to
 $\psi$ is $\beta$-\textit{compatible} with $\phi_\text{box}$ on  $X_\text{box}$,
 defined in \cite{NesNem:94}. The following assumptions will be valid throughout this section:
\begin{assumption}
\label{ass} We consider a given compact convex set $X$ with nonempty
interior and $\phi$ an associated $N_\phi$-self-concordant barrier
for $X$ (whenever $X = X_\text{box}$ we consider $\phi=
\phi_\text{box}$). Given a function $f~\in~{\cal
C}^3(\text{int}(X))$, we also assume that it satisfies one of the
properties (i)--(iii) of Lemma \ref{l_scsc}, i.e. $f$ is either
linear or convex quadratic or $M_f$-self-concordant or $X$ is a box
and $f$ satisfies condition \eqref{l_box}. Let $A \in \rset^{p
\times n}, p<n$, and $B  \in \rset^{m \times n}$ be  so that the
matrix $\bmtrx{c} A \\ B \emtrx$ has full row rank and the set $\big
\{ \{x \in \rset^n: Ax = a \} \cap \text{int}(X) \big  \} \not =
\emptyset$.
\end{assumption}

We analyze  the following prototype minimization problem:
\begin{align}
\label{prototype} \min_{x} \{ f(x) + t \phi(x) + \langle \ld , B x
\rangle: \; x \in \text{int}(X), A x =a \}.
\end{align}

 Let us define the dual convex function:
\[d(t,\ld) := \max_x \{ -f(x) - t \phi(x) - \langle \ld , B x \rangle:
\; x \in \text{int}(X), A x =a \}. \]
 Boundedness of the set $X$ and self-concordance
property of the function $f + t \phi$ (which follow from  the
assumptions mentioned above) guarantee existence and uniqueness  of
the maximizer $x(t,\ld)$ of \eqref{prototype}. Therefore, we can
consistently define the maximizer $x(t,\ld)$ of \eqref{prototype}
and the dual convex function $d(t,\ld)$   for every $t>0$ and $\ld
\in \rset^m$.

In the following four lemmas we derive the main properties of the
family of augmented dual functions $\{d(t,\cdot)\}_{t>0}$. We start
with a linear algebra result:
\begin{lemma}
\label{l0_thmain} Let $A \in \rset^{p \times n}, p<n$, and $B  \in
\rset^{m \times n}$ be two matrices and $U$ be the matrix  whose
columns  form a basis  of the null space  of $A$. Then the matrix
$\bmtrx{c} A \\ B \emtrx$ has full row rank if and only if $BU$ and
$A$ have full row rank.
\end{lemma}
\proof  Assume that $\bmtrx{c} A \\ B \emtrx$ has full row rank.
Then $A$ has full row rank.  It  remains to  show that $BU$ has full
row rank. Assume that this is not the case then  there exists a
vector $x \in \rset^{m}, x \not = 0$ such that $x^T B U =0$.  Since
the columns of $U$  span the null space of $A$ which is orthogonal
on the image space of $A^T$, it follows that $x^TB$ belongs to the
image space of $A^T$, i.e. there exists some $y \in \rset^p$ such
that $x^T B = y^T A$. But from the fact that $\bmtrx{c} A
\\ B \emtrx$ has full row rank we  must have  $x =0$ and $y=0$ which
contradicts
our assumption on $x$. \\
If $BU$ and $A$ have full row rank, it follows immediately that $B$
must have full row rank. Assume that   $\bmtrx{c} A \\ B \emtrx$
does not have full row rank. Since $A$ has full row rank, then there
exist some $y \in \rset^{m}$ and $x \in \rset^{p}, x \not = 0$, such
that \[ y^T A + x^T B = 0. \] It follows also that \[ y^T A U + x^T
B U= 0,  \quad \text{i.e.} \quad x^T B U = 0\]  and thus $x=0$ which
is a contradiction. \qed

\begin{lemma}
\label{l1_thmain} If Assumption \ref{ass} holds, then  for any $t>0$
the function $d(t,\cdot)$ is $M_t$-self-concordant, where $M_t$ is
either $2/\sqrt{t}$ or $\max\{ M_f, 2/\sqrt{t}\}$ or $2(1+
\beta)/\sqrt{t}$.
\end{lemma}
\proof Since $f$ is assumed to be either  linear or convex quadratic
or $M_f$-self-concordant or $X$ is a box and $f$ satisfies condition
\eqref{l_box} it follows from Lemma \ref{l_scsc} that $f+  t \phi$
is also $M_t$-self concordant (where $M_t$ is either $ 2/\sqrt{t}$
or $\max\{ M_f, 2/\sqrt{t}\}$ or $2(1+ \beta)/\sqrt{t}$,
respectively) and with positive definite Hessian (according to our
assumptions and Proposition \ref{l_scp}). Moreover, $f+t \phi$ is
strongly self-concordant since $\phi$ is a barrier function for $X$.
Since $A$ has full row rank and $p<n$, then there exists some
vectors $u_i, i=1 \cdots n-p$, that form a basis of the null space
of this matrix. Let $U$ be the matrix having as columns the vectors
$u_i$ and $x_0$ a particular solution of $A x= a$. Then, for a fixed
$t$, the feasible set of \eqref{prototype} can be described as

\[Q=\{y \in \rset^k: x_0 + U y \in \text{int} (X) \}, \]
which is an open convex set. Using that self-concordance is  affine
invariant it follows that the functions $\bar f(y) = f(x_0 + U y)$,
$\bar \phi(y) = \phi(x_0 + U y)$ have the same properties as the
functions $f$, $\phi$,  respectively, that $\bar f + t \bar \phi$ is
also $M_t$-self concordant and that
\[ d(t,\ld) = \max_{y \in Q} [ - \bar f(y) - t \bar \phi(y) -  \langle \ld
, B (x_0 + U y) \rangle]. \]

From our assumptions and Proposition \ref{l_scp} it follows that the
Hessian of $\phi$ and $\bar \phi$ are positive definite. Since $f$
is convex it follows that the Hessian of $\bar f+ t \bar \phi$ is
also positive definite and thus invertible. Let \[ \bar F (t,\ld) =
\max_{y \in Q} [ \langle \ld , y \rangle - \bar f(y) - t \bar
\phi(y)] \] be the Legendre transformation of $\bar f  + t \bar \phi
$. In view of known properties of the Legendre transformation, it
follows that if $\bar f + t \bar \phi$ is convex on $X$ from ${\cal
C}^3$ such that its Hessian is positive definite, then  $\bar
F(t,\cdot)$ has the same properties on its domain $\{ \ld \in \rset
^m: \langle \ld , y \rangle - \bar f(y) - t \bar \phi(y) \  \text{
bounded above on} \ Q \}$. Moreover,  from Theorem 2.4.1 in
\cite{NesNem:94} it follows that $\bar F(t,\cdot)$ is also
$M_t$-self-concordant on its domain. Note that \[ d(t,\ld) =
 \langle \ld , - B x_0 \rangle + \bar F(t,-(BU)^T \ld). \]
  Since
 $\bmtrx{c} A \\ B \emtrx$ has full row rank,  then from Lemma \ref{l0_thmain}
 $BU$  has full row rank. Moreover, since $\nabla^2 \bar F(t,\cdot)$  is positive
 definite and  \[ \nabla^2 d(t,\ld) = BU \nabla^2\bar F \big(t,-(BU)^T \ld \big ) (BU)^T, \]
  it follows that $\nabla^2 d(t,\cdot)$ is positive definite on its
domain
 \[ X_{d(t,\cdot)}:= \{ \ld \in \rset ^m: - \bar f(y) - t \bar \phi(y)
- \langle \ld , B (x_0 + U y) \rangle \ \text{ bounded above on} \ Q
\}. \]
 Moreover, since self-concordance is  affine invariant  it
follows that $d(t,\cdot)$ is also
  $M_t$-self-concordant on the domain $X_{d(t,\cdot)}$. \qed

\begin{lemma}
\label{l2_thmain}  Under Assumption \ref{ass} the  inequality $|
\langle \nabla d^\prime (t,\ld), h \rangle | \leq (2 \xi_t/\alpha_t)
\sqrt{h^T \nabla^2 d(t,\ld) h } $ \ holds true for each $t>0$ and
$\ld, h \in \rset^m$, where $\xi_t=(M_t/2) \sqrt{N_\phi / t}$  and
$\alpha_t = M_t$.
\end{lemma}
\proof From Lemma \ref{l1_thmain} we know that $d(t,\cdot)$ is
${\cal C}^3$ with positive definite Hessian. By virtue of the
barrier $\phi$ for the set $X$ the optimal solution  $x(t,\ld)$ of
\eqref{prototype} satisfies $x(t,\ld) \in \text{int}(X)$ and so the
first-order optimality conditions  for optimization problem
\eqref{prototype} are: there exists $\nu(t,\ld) \in \rset^p$ such
that
\begin{align}
\label{kkt}
 \nabla f(x(t,\ld)) + t \nabla
\phi(x(t,\ld)) + B^T \ld + A^T \nu(t,\ld) = 0
 \; \text{and} \; A x(t,\ld) = a.
\end{align}
 First we determine the formula for the Hessian. It follows
immediately from \eqref{kkt} that
\[ \nabla d(t,\ld) = - B x(t,\ld) \quad \text{and} \quad
 \nabla^2 d(t,\ld) = - B \nabla x(t,\ld).   \]
 Let us introduce the following notation:
\[ H(t,\ld):= \nabla^2 f(x(t,\ld)) + t \nabla^2
\phi(x(t,\ld)).  \]
 For simplicity,   we drop
the dependence of all the functions on $x(t,\ld)$ and  $(t,\ld)$.
Differentiating    \eqref{kkt} with respect to $\ld$ we arrive at
the following system in $\nabla x$ and $\nabla \nu$:
\[ \bmtrx{cc}
      \nabla^2 f + t \nabla^2 \phi     &   A^T   \\
            A                         &    0
   \emtrx
   \bmtrx{c}
      \nabla x   \\
       \nabla \nu
   \emtrx =
   \bmtrx{c}
       - B^T \\
         0
   \emtrx.
\]
Since $H$ is positive definite and according to our assumption $A$
is full row rank, it follows that the system matrix is invertible.
Using the well-known formula for inversion of partitioned matrices
we find that:
\begin{align}
\label{form_hessian} \nabla^2 d = B [H^{-1} - H^{-1} A^T (A
H^{-1}A^T )^{-1} A H^{-1}]B^T .
\end{align}
 Differentiating the first part of \eqref{kkt} with respect
to $t$ and using the same procedure as before we arrive at a similar
system as above in the unknowns $x^\prime$ and $\nu^\prime$. We find
that
\[ x^\prime = - [H^{-1} - H^{-1} A^T (A H^{-1}A^T )^{-1} A H^{-1}] \nabla \phi \quad
\text{and} \quad  \nabla d^\prime  = - B x^\prime.\]
 We also introduce the following notation:
$F:= H^{-1} A^T (A H^{-1}A^T )^{-1} A H^{-1}$ and $G: = H^{-1} - F$,
which are  positive semidefinite. Using a similar reasoning as in
\cite{Zha:05} and Cauchy-Schwarz inequality we obtain:
\begin{align*}
| \langle \nabla d^\prime, h \rangle | =  | h^T B G \nabla \phi|
\leq \sqrt{h^T B G B^T h} \sqrt{\nabla \phi^T G  \nabla \phi} =
\sqrt{h^T (\nabla^2d) h} \sqrt{\nabla \phi^T G  \nabla \phi}.
\end{align*}
 Since \[ G= H^{-1} - F \preceq H^{-1} = (\nabla^2 f + t
\nabla^2 \phi)^{-1} \ \preceq 1/t (\nabla^2 \phi)^{-1} \] and using
\eqref{d_scb1} it follows that \[ | \langle \nabla d^\prime, h
\rangle | \leq  \sqrt{h^T (\nabla^2d) h} \sqrt{N_\phi / t}. \] \qed

\begin{lemma}
\label{l3_thmain} Under Assumption \ref{ass} the  inequality  $|
\langle  \nabla^2 d^\prime (t,\ld) h, h \rangle | \leq 2 \eta_t \
h^T \nabla^2 d(t,\ld) h$ \ holds true for each $t>0$ and $\ld, h \in
\rset^m$, where $\eta_t=(M_t/2) \sqrt{N_\phi/t} + (1/2t)$.
\end{lemma}
\proof  We recall that $H(t,\ld) = \nabla^2 f(x(t,\ld)) + t \nabla^2
\phi(x(t,\ld))$. Therefore
\[ \begin{array}{c}
  h^T H^\prime(t,\ld) h = (\nabla^3
f(x(t,\ld)) + t \nabla^3 \phi(x(t,\ld)))[x^\prime(t,\ld), h, h] +
h^T \nabla^2 \phi(x(t,\ld)) h, \\
 h^T(H^{-1}(t,\ld))^\prime h = -h^T H^{-1}(t,\ld)  H^\prime(t,\ld)
H^{-1}(t,\ld) h.
\end{array}
\]
We again drop the dependence on $(t,\ld)$ and after some
straightforward algebra computations we arrive at the following
expression:
\[ \langle  \nabla^2 d^\prime  h, h \rangle = - h^T B (H^{-1} - F) H^\prime (H^{-1} - F) B^T h. \]
Let us denote with \[ u := (H^{-1} - F) B^T
h. \]
Taking into account the expression of $H^\prime$ derived above
we obtain: \[ |\langle  \nabla^2 d^\prime  h, h \rangle | = |u^T
H^\prime u| = |(\nabla^3 f + t \nabla^3 \phi)[x^\prime, u, u] + u^T
\nabla^2 \phi u|. \]
 Using the self-concordance property
\eqref{d_sce} for $f + t \phi$ we obtain that:
\begin{align*} |(\nabla^3 f + t \nabla^3
\phi)[x^\prime, u, u] | & \leq  M_t u^T(\nabla^2 f + t \nabla^2
\phi)u \sqrt{(x^\prime)^T(\nabla^2 f + t \nabla^2 \phi)x^\prime }
\\
&= M_t u^T H u \sqrt{(x^\prime)^T H x^\prime}.
\end{align*}

 Moreover, since $f$ is convex, $\nabla^2f$
is positive semidefinite
 and thus:  \[ u^T \nabla^2 \phi u \leq (1/t) \ u^T H u. \]
 Combining the last two inequalities we
obtain:
\begin{align}
\label{l_ineq3} |\langle \nabla^2 d^\prime  h, h \rangle |  \leq M_t
u^T H u \sqrt{(x^\prime)^T H x^\prime} +  (1/t)  u^T H u .
\end{align}

With some algebra we can check that the following identity holds: $F
H (H^{-1}~-~F)~=~0$. Based on this identity we can compute $u^T H u$
and $(x^\prime)^T H x^\prime$. Indeed,
\begin{align*}
& u^T H u  = h^T B (H^{-1} - F) H (H^{-1} - F)B^T h  \\
& = h^T B (H^{-1} - F) B^T h - h^T B F H (H^{-1} - F)B^T h =  h^T B
(H^{-1} - F) B^T h = h^T \nabla^2 d \ h.
\end{align*}
 Similarly, using \eqref{d_scb1} we obtain
\begin{align*}
& (x^\prime)^T H x^\prime = \nabla \phi^T   (H^{-1} - F) \nabla \phi
\leq \nabla \phi^T   H^{-1}  \nabla \phi \leq (1/t)   \nabla \phi^T
(\nabla^2 \phi)^{-1}  \nabla \phi \leq N_\phi/t.
\end{align*}
The inequality from  lemma  follows then by replacing the last two
relations in \eqref{l_ineq3}. \qed

The main result of this section  is summarized in the next theorem.
\begin{theorem}
\label{th_main}
 Under the   Assumption \ref{ass},
$\{d(t,\ld)\}_{t>0}$ is a strongly self-concordant family in the
sense of Definition\footnote{Note that according to Definition 3.1.1
in \cite{NesNem:94}  $\gamma_t=1$ and $ \mu_t=1$  for our case.}
3.1.1 in \cite{NesNem:94} with parameters $\alpha_t=M_t, \xi_t=
(M_t/2) \sqrt{N_\phi/t}$ and $\eta_t=(M_t/2) \sqrt{N_\phi/t} +
(1/2t)$, where $M_t$ is defined in Lemma \ref{l1_thmain}.
\end{theorem}
\proof  Basically, from Definition 3.1.1 in \cite{NesNem:94} we must
check three properties: self-concordance of  $d(t,\ld)$ (Lemma
\ref{l1_thmain}) and that the first and second order derivative of
$d(t,\cdot)$ vary with $t$ at a rate proportional to the derivative
itself (Lemmas \ref{l2_thmain} and \ref{l3_thmain}). In conclusion,
the Lemmas \ref{l1_thmain}--\ref{l3_thmain} prove our theorem. \qed


It is known \cite{NesNem:94} that
self-concordant families of functions can be minimized by
path-following methods in polynomial time. Therefore, this type of
family of augmented dual functions $\{d(t,\cdot)\}_{t>0}$ plays an
important role in the algorithm of the next section.


\section{Parallel Implementation of an Interior-Point Based Decomposition Method}
\label{s_pflda}

In this section we develop an interior-point Lagrangian
decomposition method for the separable convex  problem given by
\eqref{scp_f}--\eqref{scp_B}. Our previous Theorem \ref{th_main} is
the major contribution of our paper since it allows us to
effectively utilize the  Newton method for tracing the trajectory of
optimizers of the self-concordant family of augmented dual functions
\eqref{dip}.

\subsection{Interior-Point Lagrangian Algorithm}

The following assumptions for optimization problem
\eqref{scp_f}--\eqref{scp_B} will be valid in this section:
\begin{assumption}
\label{ass1} (i) The sets $X_i$ are compact convex sets with
nonempty interior and $\phi_{X_i}$ are $N_i$-self-concordant
 barriers for $X_i$.\\
(ii) Each function $f_i$ is either linear or convex quadratic or
$M_{f_i}$-self-concordant or $X_i$ is a box and $f_i$ satisfies
condition \eqref{l_box}. \\
(iii) The   \textit{block-angular} matrix $\bmtrx{c}
      D_A    \\ B \emtrx$   has full row rank  and the set
   $\big \{\{x \in \rset^n: A_ix_i = a_i, \ B x =
b \} \cap \text{int}(X)\big \} \not = \emptyset$.
\end{assumption}
Note that  boundedness of the set $X_i$ can be relaxed to $X_i$ does
not contain straight lines and  the set of  optimal solutions
 to problem  \eqref{scp_f}--\eqref{scp_B} is bounded. Note also that the rank
 assumption (iii)  is not restrictive since we can eliminate the redundant equalities
 (see also Lemma \ref{l0_thmain} for other less restrictive conditions).  The constraint qualification condition
 from Assumption \ref{ass1}
(iii) guarantees that strong duality holds for problem
\eqref{scp_f}--\eqref{scp_B} and thus there exists a primal-dual
optimal solution $(x^*, \ld^*)$.

Let us  introduce the dual function:
\begin{align*}
  d(t,\ld) &=  \max_{x} \{
-L_t^{\text{sc}}(x,\ld) : x_i \in \text{int}(X_i),  \; A_i x_i = a_i
\; \forall i=1 \cdots N \}  \\
& = \langle \ld, b \rangle + \sum_{i=1}^N \max_{x_i} \{-f_i(x_i) - t
\phi_{X_i} (x_i) - \langle \ld,B_i x_i \rangle  : x_i \in
\text{int}(X_i),  \; A_i x_i = a_i\} \\
& = \langle \ld, b \rangle + \sum_{i=1}^N d_i(t,\ld).
\end{align*}
 Note that the  function $d(t,\cdot)$ can
be computed in \textit{parallel} by decomposing the original large
optimization problem \eqref{scp_f}--\eqref{scp_B} into $N$
independent small convex subproblems.
\begin{lemma}
\label{l_sscfd} (i) The family $\{ d_i(t,\cdot) \}_{t>0}$ is
strongly self-concordant with the parameters $\alpha_i(t)=M_i(t),
\xi_i(t)= (M_i(t)/2) \sqrt{N_i/t}$ and $\eta_i(t)=(M_i(t)/2)
\sqrt{N_i/t} + (1/2t)$, where $M_i(t)$ is  either $2/\sqrt{t}$ or
$\max\{ M_{f_i},
2/\sqrt{t}\}$ or $2(1+ \beta)/\sqrt{t}$ for all $i=1 \cdots N$.\\
(ii) The family $\{ d(t,\cdot) \}_{t>0}$ is strongly self-concordant
with parameters $\alpha(t) = \alpha/\sqrt{t}$,  $\xi(t) = \xi/t$ and
$\eta(t) = \eta/ t$, for some fixed positive constants $\alpha, \xi$
and $\eta$ depending on $(N_i, M_{f_i},\beta)$.
\end{lemma}
\proof (i) is a straightforward consequence of  Assumption
\ref{ass1} and Theorem \ref{th_main}.

(ii) Note that $d(t,\ld) = \langle \ld, a \rangle + \sum_{i=1}^N
d_i(t,\ld)$.  From Proposition 3.1.1 in \cite{NesNem:94} we have
that the sum of strongly self-concordant family of functions is also
strongly self-concordant family with the parameters: $\alpha(t) \geq
\max_i \{ \alpha_i(t)\}$ is a positive continuously differentiable
function on $\rset_+$, $ \xi(t) = \alpha(t) \max_i \{ 2 \xi_i(t)/
\alpha_i(t) \}$ and $\eta(t) = \max_i \{ \eta_i(t)\}$. Since
$\alpha_i(t)$ is  either $2/\sqrt{t}$ or  $\max\{ M_{f_i},
2/\sqrt{t}\}$ or $2(1+ \beta)/\sqrt{t}$ for all $i=1 \cdots N$, it
follows that we can always choose $\alpha(t) = \alpha/\sqrt{t}$,
where $\alpha = 2$ or $\alpha = 2(1+\beta)$. Similarly, we can show
that there exists positive constants $\xi$ and $\eta$ depending on
$N_i, M_{f_i}$ and $\beta$ such that $\xi(t) = \xi/t$ and $\eta(t) =
\eta/ t$. \qed

 From Assumption  \ref{ass1} and the discussion from
previous section, it follows that the optimizer of each maximization
is unique and denoted by
\begin{align}
\label{eq_exact} x_i(t,\ld) := \arg \max_{x_i} \{-f_i(x_i) - t
\phi_{X_i} (x_i) - \langle \ld,B_i x_i \rangle  : x_i \in
\text{int}(X_i),  \; A_i x_i = a_i\}
\end{align}
 and $x(t,\ld):=[x_1(t,\ld)^T \cdots x_N(t,\ld)^T]^T$. It is clear that the
augmented dual function $d^{\text{sc}}(t,\ld) = - d(t,\ld)$ and let
$\ld(t) := \arg \max_{\ld \in \rset^m} d^{\text{sc}}(t,\ld)$,  or
equivalently
\[ \ld(t)= \arg \min_{\ld \in \rset^m} d(t,\ld). \]

 From Assumption \ref{ass1} and the proof of
Lemma \ref{l1_thmain} it follows that the Hessian $\nabla^2
d(t,\ld)$ is positive definite for all $t>0$ and $\ld \in \rset^m$.
Hence, the dual function $d(t,\cdot)$ is strictly convex and thus
$\ld(t)$ is unique. Therefore, we can consistently define the set
$\{ (x(t,\ld(t)),\ld(t)): t>0 \}$, called the \textit{central path}.
Let us introduce the $N_\phi$-self-concordant barrier function
$\phi_X(x):=\sum_{i=1}^N \phi_{X_i}(x_i)$ for the set $X$, where
$N_\phi= \sum_{i=1}^N N_i$.
\begin{lemma}
\label{l_convct} The central path $\{ (x(t,\ld(t)), \ld(t)): t>0 \}$
converges to the  optimal solution $(x^*, \ld^*)$ as $t \to 0$ and
$\{x(t,\ld(t)): t>0 \}$  is  feasible for the  problem
\eqref{scp_f}--\eqref{scp_B}.
\end{lemma}
\proof Let $x(t):= \arg \min_x \{ f(x) + t \phi_X(x): Bx = b,  x_i
\in \text{int}(X_i),   A_i x_i = a_i \; \forall i \}$, then it is
 known that $x(t) \to x^*$ as $t \to 0$. It is easy to see that the
Hessian of $f + t \phi_X$ is positive definite and thus $f + t
\phi_X$ is strictly convex and $x(t)$ is unique. From Assumption
\ref{ass1} it also follows that strong duality holds for this
barrier function problem and therefore
\begin{align*}
&\min_x \{ f(x) + t \phi_X(x): Bx = b,  x_i \in \text{int}(X_i), A_i
x_i = a_i \; \forall i\}= \\
& \max_\ld \min_x \{ f(x) + t \phi_X(x) + \langle \ld, Bx - b
\rangle : x_i \in \text{int}(X_i), A_i x_i = a_i \; \forall i\} = \\
&  \min_x \{ f(x) + t \phi_X(x) + \langle \ld(t), Bx - b \rangle :
x_i \in \text{int}(X_i), A_i x_i = a_i \; \forall i\}.
\end{align*}
 In conclusion, $x(t) = x(t,\ld(t))$ and thus $x(t, \ld(t))
\to x^*$ as $t \to 0$. As a consequence it follows that $x(t,
\ld(t))$ is feasible for the original problem, i.e. $B
x(t,\ld(t))~=~b$, $A_i x_i(t,\ld(t)) = a_i$ and $x_i(t,\ld(t)) \in
\text{int}(X_i)$. It is also clear that $\ld(t) \to \ld^*$ as $t \to
0$.\qed

The next theorem  describes the behavior  of the central path:
\begin{theorem}
\label{th_convrates} For  $x(t) = x(t,\ld(t))$ the following bound
holds for the central path: given any $0 < \tau < t$ then,
\[  f(x(t)) - f(x(\tau)) \leq  N_\phi (t-\tau). \]
\end{theorem}
\proof  For any $s>0$, $x(s)=[x_1(s)^T \cdots x_N(s)^T]^T$ satisfies
the following optimality conditions (see \eqref{kkt} and Lemma
\ref{l_convct}): there exists  $\nu(s) \in \rset^{\sum_{i=1}^N m_i}$
such that
\[ \nabla f(x(s)) + s \nabla \phi_X(x(s)) + B^T \ld(t) +
D_A^T \nu(s) = 0, \;\ B x(s) = b \;\ \text{and} \;\ A_i x_i(s) =
a_i.
\]

 It follows immediately that $\langle
\nabla f(x(s)), x^\prime(s) \rangle = - s \langle \nabla
\phi_X(x(s)), x^\prime(s) \rangle$. Since $0< \tau < t$, then there
exists $\tau \leq s \leq t$ such that
\begin{align*}
f(x(t)) -f(x(\tau)) = (t-\tau) \langle \nabla f(x(s)), x^\prime(s)
\rangle = -s(t-\tau) \langle \nabla \phi_X(x(s)), x^\prime(s)
\rangle.
\end{align*}

From \eqref{d_scb} we have that
\[ - \langle \nabla \phi_X(x(s)), x^\prime(s) \rangle \leq
\big(N_\phi \ x^\prime(s)^T \nabla^2 \phi_X(x(s)) x^\prime(s)
\big)^{1/2}.\]   Using a similar reasoning as in Lemma
\ref{l2_thmain} we have:
\[ x^\prime(s) = - [H^{-1}(s) - H^{-1}(s) D_A^T  (D_A H^{-1}(s)D_A^T )^{-1} D_A H^{-1}(s)]
 \nabla \phi_X(x(s)),  \]
where  we denote with $H(s) =\nabla^2 f(x(s)) + s \nabla^2
\phi_X(x(s))$. Using \eqref{d_scb1}, the expression for
$x^\prime(s)$ and since $0 \prec \nabla^2 \phi_X(x(s)) \preceq 1/s
H(s)$ and $H^{-1}(s) \preceq 1/s \big(\nabla^2
\phi_X(x(s))\big)^{-1}$ we obtain:
\begin{align*}
 x^\prime(s)^T \nabla^2 \phi_X(x(s)) x^\prime(s) & \leq  (1/s)
\nabla \phi_X(x(s))^T  H^{-1}(s)
  \nabla \phi_X(x(s))   \\
  & \leq (1/s^2)   \nabla \phi_X(x(s))^T  \big(\nabla^2 \phi_X(x(s))\big)^{-1}
  \nabla \phi_X(x(s)) \leq N_\phi / s^2.
\end{align*}

 It follows immediately  that $ f(x(t))
-f(x(\tau)) \leq N_\phi (t-\tau)$. \qed

A simple consequence of Theorem \ref{th_convrates} is that the
following bounds on the approximation of the optimal value function
$f^*$ hold:
\[ 0 \leq f(x(t)) - f^* \leq  t N_\phi. \]
 Indeed, from Lemma \ref{l_convct} we know
that $\{x(t,\ld(t)): t>0 \}$ is feasible for the original problem
\eqref{scp_f}--\eqref{scp_B}. Since $x(t)= x(t,\ld(t))$, it follows
that $f(x(t)) \geq f^*$. It remains to show the upper bound.
However, taking the limit as $\tau \to 0$ in Theorem
 \ref{th_convrates} and using Lemma
\ref{l_convct} we obtain also the upper bound. This  upper bound
 gives us a stopping criterion in the algorithm that we derive below:
 if $\epsilon$ is the required  accuracy  for the
approximation of $f^*$, then for any $t_f \leq \epsilon/N_\phi$ we
have  that $x(t_f)$ is an $\epsilon$-approximation of the optimum,
i.e. $x(t_f)$ is feasible for problem \eqref{scp_f}--\eqref{scp_B}
and $f(x(t_f)) - f(x^*) \leq \epsilon$. Although $\ld(t)$ is the
minimizer of the dual function $d(t,\cdot)$ over $\rset^m$, so
various unconstrained minimization techniques (e.g. Newton,
quasi-Newton and conjugate gradient methods) can be used to
approximate $\ld(t)$, our goal is to trace the central path
$\{(x(t,\ld(t)),\ld(t)):  t>0 \}$ utilizing Newton method for the
self-concordant family $\{d(t,\cdot)\}_{t>0}$.

It is easy to see that the gradient of the self-concordant function
$d(t,\cdot)$ is given by
\[\nabla d(t,\ld) = b + \sum_{i=1}^N \nabla d_i(t,\ld)   =
b - \sum_{i=1}^N B_i x_i(t,\ld) = b - B x(t,\ld). \]  For every
$(t,\ld)$ let us define the positive definite matrix
\[ H_i(t,\ld) := \nabla^2 f_i(x_i(t,\ld)) + t \nabla^2 \phi_{X_i}(x_i(t,\ld)). \]

 The Hessian of  function
$d_i(t,\cdot)$ is positive definite and from \eqref{form_hessian} it
has the  form
\[ \nabla^2 d_i(t,\ld) =  B_i[H_i(t,\ld)^{-1} -
 H_i(t,\ld)^{-1} A_i^T \big(A_i H_i(t,\ld)^{-1}
 A_i^T \big)^{-1} A_i H_i(t,\ld)^{-1}  ]B_i^T.  \]

In conclusion, the Hessian of  dual function $d(t,\cdot)$ is also
positive definite and given by:
\[ \nabla^2 d(t,\ld) = \sum_{i=1}^N \nabla^2 d_i(t,\ld).  \]
 Denote the Newton direction
 associated to  self-concordant function $d(t, \cdot)$ at $\ld$ with

\[ \Delta \ld (t,\ld) := - \big (\nabla^2 d(t,\ld) \big )^{-1} \nabla d(t,\ld).  \]

 For every $t > 0$, we define the Newton
decrement of the function $d(t,\cdot)$ at $\ld$ as:

\[  \delta(t,\ld) := \alpha(t)/2 \sqrt{ \nabla d(t,\ld)^T
 \big ( \nabla^2 d(t,\ld) \big )^{-1}  \nabla d(t,\ld) }. \]
 Note that $\delta(t,\hat \ld) = 0$ if and only if $\hat
\ld =\ld(t)$ (recall that $ \ld(t) = \arg \min_{\ld \in \rset^m}
d(t,\ld)$).

\begin{algorithm}(Initialization of  Path-Following Algorithm)\\
\label{alg_ip1}
Step 0. \textbf{input} $t_0>0$, $\ld_0  \in \rset^m$, $\epsilon_V>0$
and
$r=0$\\
Step 1. compute $x_i^r \!=\! x_i(t_0, \ld_r) \ \forall i$, $\delta_r
\!=\!
\delta(t_0,\ld_r)$; if $\delta_r \leq \epsilon_V$, $r_f = r$ and go to Step 3\\
Step 2. determine a step size $\sigma$ and compute  Newton iterate:
$\ld_{r+1}= \ld_r + \sigma \Delta \ld(t_0,\ld_r)$;

\hspace*{0.3cm} replace $r$ by $r+1$ and go to Step 1 \\
Step 3. \textbf{output} $(t^0,\ld^0) = (t_0,\ld_{r_f})$.
\end{algorithm}
Note that Algorithm \ref{alg_ip1} approximates the optimal Lagrange
multiplier $\ld(t_0)$ of the dual function $d(t_0,\cdot)$, i.e. the
sequence $(t_0,\ld_r)$ moves into the neighborhood $V(t,\epsilon_V)
=\{(t,\ld): \delta(t,\ld) \leq \epsilon_V \}$ of the trajectory
$\{(t,\ld(t)): t >0\}$.

\begin{algorithm}(Path-Following Algorithm)\\
\label{alg_ip2} Step 0. \textbf{input}: $(t^0,\ld^0)$ satisfying
$\delta(t^0,\ld^0) \leq \epsilon_V$ , $k=0$, $0 < \tau < 1$ and
$\epsilon>0$ \\
Step 1. if $t^k N_\phi \leq \epsilon$, then $k_f=k$ and go to Step 5\\
Step 2. (\underline{outer iteration}) let $t^{k+1} = \tau t^k$ and
go to
inner iteration (Step 3)\\
Step 3. (\underline{inner iteration}) initialize $\ld = \ld^k$,
$t=t^{k+1}$ and $\delta =
\delta(t^{k+1},\ld^k)$\\
while $\delta > \epsilon_V$ do

Step 3.1 compute $x_i = x_i(t, \ld) \; \forall i$, determine a step
size $\sigma$ and compute

\hspace*{1.5cm} $\ld^+ = \ld + \sigma \Delta \ld(t,\ld)$

Step 3.2 compute $\delta^+ = \delta(t,\ld^+)$ and update $\ld=\ld^+$ and $\delta = \delta^+$\\
Step 4. $\ld^{k+1} = \ld$ and $x_i^{k+1} = x_i$; replace $k$ by
$k+1$ and
go to Step 1\\
Step 5. \textbf{output}: $(x_1^{k_f}, \cdots, x_N^{k_f},
\ld^{k_f})$.
\end{algorithm}
In Algorithm \ref{alg_ip2}  we trace numerically the trajectory
$\{(t,\ld(t)): t >0\}$ from a given initial point $(t^0,\ld^0)$
close to this trajectory.  The sequence $\{(x_1^{k}, \cdots,
x_N^{k}, \ld^k )\}_{k>0}$ lies in a neighborhood of the central path
and each limit point of this sequence is primal-dual optimal.
Indeed, since $t^{k+1} = \tau t^k$  with $\tau < 1$, it follows that
$\lim_{k \to \infty} t^k = 0$ and using Theorem \ref{th_convrates}
the convergence of the sequence $x^k =[(x_1^k)^T \cdots
(x_N^k)^T]^T$ to $x^*$ is obvious.

The step size $\sigma$ in the previous algorithms is defined by some
line search rule. There are many strategies for choosing $\tau$.
Usually, $\tau$ can be chosen independent  of the problem (long step
methods), e.g. $\tau = 0.5$, or depends on the problem (short step
methods). The choice for $\tau$ is crucial for the performance of
the algorithm. An example is that in practice long step
interior-point algorithms are more efficient than short step
interior-point algorithms. However, short step type algorithms have
better worst-case complexity iteration bounds than long step
algorithms. In the sequel we derive a theoretical strategy to update
the barrier parameter $\tau$ which follows from the theory described
in \cite{NesNem:94} and consequently we obtain complexity bounds for
short step updates. Complexity iteration bounds for long step
updates can also be derived using the same theory (see Section 3.2.6
in \cite{NesNem:94}). The next lemma estimates the reduction of the
 dual function at each  iteration.
\begin{lemma}
\label{l_anal} For any $t>0$ and $\ld  \in \rset^m$, let $\Delta \ld
= \Delta \ld (t,\ld)$ be the Newton direction as defined above. Let
also $\delta = \delta (t,\ld)$ be the Newton decrement and $\delta_*
= 2 - \sqrt{3}$. \\
(i) If $\delta > \delta_*$, then defining the step length $\sigma =
1/(1+\delta)$ and the Newton iterate $\ld^+ = \ld + \sigma \Delta
\ld$ we have the following decrease in the objective function
$d(t,\cdot)$
\[ d(t,\ld^+) - d(t, \ld) \leq -(4 t/\alpha^2)
(\delta - \log (1+\delta)).  \]

 (ii) If $\delta \leq \delta_*$, then defining the
Newton iterate $\ld^+ = \ld + \Delta \ld$ we have
\[ \delta(t,\ld^+) \leq \delta^2/(1-\delta)^2 \leq \delta/2, \quad
d(t,\ld) - d(t, \ld(t)) \leq  (16 t/\alpha^2) \delta.
\]

 (iii)  If $\delta \leq \delta_*/2$, then defining
$t^+ =\frac{2 c}{2 c +1}  t$, where $c = 1/4 + 2 \xi /\delta_* +
\eta$,  we have
\[ \delta(t^+,\ld) \leq \delta_*. \]
\end{lemma}
 \proof (i) and (ii) follow from Theorem 2.2.3 in
\cite{NesNem:94} and Lemma \ref{l_sscfd} from above.

(iii) is based on the result of Theorem 3.1.1 in \cite{NesNem:94}.
In order to apply this theorem,  we first write the metric defined
by (3.1.4) in \cite{NesNem:94} for our problem: given $0 < t^+ < t$
and using Lemma \ref{l_sscfd} we obtain
\[ \rho_{\delta_*/2}(t, t^+) = (1/4 + 2 \xi /\delta_* + \eta) \log(t/t^+).  \]

 Since $\delta \leq \delta_*/2 <  \delta_*$ and since
for $t^+ = \frac{2 c}{2 c +1} t$, where $c$ is defined above, one
can verify that $\rho_{\delta_*/2}(t, t^+) = c \log(1 + 1/2c) \leq
1/2 \leq 1- \delta/\delta_*$, i.e. $t^+$ satisfies the condition
(3.1.5) of Theorem 3.1.1 in \cite{NesNem:94}, it follows that
$\delta(t^+,\ld) \leq \delta_*$. \qed

Define the following step  size: $\sigma(\delta) = 1/(1+\delta)$ if
$\delta  > \delta_*$ and  $\sigma(\delta) = 1$ if $\delta \leq
\delta_*$. With Algorithm \ref{alg_ip1} for a given $t^0$ and
$\epsilon_V= \delta_*/2$, we can find $(t^0,\ld^0)$ satisfying
$\delta(t^0,\ld^0) \leq \delta_*/2$ using the step size
$\sigma(\delta)$ (see previous lemma). Based on the analysis given
in Lemma \ref{l_anal} it follows that taking in Algorithm
\ref{alg_ip2} $\epsilon_V = \delta_*/2$ and $\tau = 2 c/ (2c +1)$,
then the inner iteration stage (step 3)  reduces to only one
iteration:

\textit{Step 3. \text{compute} $\ld^{k+1} = \ld^k + \Delta
\ld(t^{k+1},\ld^k)$.}

 However, the number of outer iterations is larger than in
the case of long step algorithms.

\subsection{Practical Implementation}
\label{s_PI}

In this section we discuss the practical implementation of our
algorithm and we give some estimates of the complexity  for it.
Among the assumptions considered until now in the paper the most
stringent one seems to be the one requiring to solve exactly the
maximization problems \eqref{eq_exact}, i.e. the exact computation
of the maximizers $x_i(t,\ld)$. Note that the gradient and the
Hessian of $d(t,\cdot)$ at $\ld$ depends on $x_i(t,\ld)$'s. When
$x_i(t,\ld)$'s are computed approximately, the expressions for the
gradient and Hessian derived in the previous section for
$d(t,\cdot)$ at $\ld$ are not the true gradient and Hessian of
$d(t,\cdot)$ at this point. In simulations we considered the
following criterion:  find $\tilde x_i(t,\ld) \in \text{int}(X_i)$
and $\tilde \nu_i(t,\ld) \in \rset^{m_i}$ such that $A_i \tilde
x_i(t,\ld) = a_i$ and the following condition holds
\[ \| \nabla f_i(\tilde x_i(t,\ld)) + t \nabla \phi_{X_i}(\tilde x_i(t,\ld)) +
B_i^T \ld + A_i^T  \tilde \nu_i(t,\ld)  \| \leq t \epsilon_x, \] for
some $\epsilon_x >0$. Note however that even when such
approximations are considered, the vector $\Delta \ld$ still
 defines a search direction in the $\ld$-space. Moreover, the cost
 of computing an extremely  accurate maximizer of \eqref{eq_exact}
 as compared to the cost of  computing a good maximizer of
 \eqref{eq_exact} is only marginally more, i.e. a few Newton steps
 at most (due to quadratic convergence of the Newton method close to the
 solution). Therefore,  it is not unreasonable to assume even exact
 computations in the proposed algorithms.

\subsubsection{Parallel Computation}

In the rest of this section we discuss the complexity of our method
and  parallel implementations for solving efficiently the Newton
direction  $\Delta \ld$. At each iteration of the algorithms  we
need to solve basically a linear system of the following form:
\begin{align}
\label{nls}
 \big( \sum_{i=1}^n G_i \big) \Delta \ld =g,
 \end{align}
  where $G_i= B_i[H_i^{-1} -
 H_i^{-1} A_i^T \big(A_i H_i^{-1}
 A_i^T \big)^{-1} A_i H_i^{-1}  ]B_i^T $, the positive definite matrix
$H_i$ denotes the Hessian of $f_i + t  \phi_{X_i}$ and some
appropriate vector $g$. In order to obtain the matrices $H_i$ we can
solve in parallel $N$ small convex optimization problems of the form
\eqref{eq_exact} by Newton method, each one of dimension $n_i$  and
with self-concordant objective function. The cost to solve each
subproblem \eqref{eq_exact} by Newton method is ${\cal O}
(n_i^3(n_\ld + \log \log 1/t \epsilon_x))$, where $n_\ld$ denotes
the number of Newton iterations before the iterates $x_i$ reaches
the quadratic convergence region (it depends on the update $\ld$)
and $t \epsilon_x$ is the required accuracy for the approximation of
\eqref{eq_exact}. Note that using the Newton method for solving
\eqref{eq_exact}   we automatically obtain also the expression for
$H_i^{-1}$  and $A_i H_i^{-1}  A_i^T$. Assuming that a Cholesky
factorization  for $A_i H_i^{-1}  A_i^T$ is used to solve the Newton
system  corresponding to the optimization subproblem
\eqref{eq_exact}, then this factorization can also be used to
compute in parallel the matrix of the linear system \eqref{nls}.
Finally, we can use a Cholesky factorization of this matrix and then
forward and  backward substitution to obtain the Newton direction
$\Delta \ld$. In conclusion, we can compute the Newton direction
$\Delta \ld$ in ${\cal O} (\sum_{i=1}^N n_i^3)$ arithmetic
operations.

Note however that in many applications the matrices $H_i$,  $A_i$
and $B_i$ are very sparse and have special structures. For example
in network optimization (see Section \ref{sec_no}  below for more
details) the $H_i$'s are diagonal matrices, $B_i$'s
 are the identity matrices and the matrices $A_i$'s are the same for
 all $i$ (see \eqref{general_dcop}), i.e. $A_i=A$. In this case the Cholesky
 factorization of $A H_i^{-1}  A^T$ can be done very
 efficiently since the sparsity pattern of those matrices  is the same in all iterations
 and coincides with the sparsity pattern of $A A^T$,
 so the \textit{analyse} phase has to be done only once, before optimization.

 For large problem instances we can also  solve the linear system \eqref{nls} approximately
 using a preconditioned  conjugate gradient algorithm. There are
 different techniques to construct a good preconditioner  and
 they  are spread across optimization literature. 
Detailed simulations for the   method proposed in this paper and
comparison of different techniques to solve the Newton system
\eqref{nls} will be given elsewhere.

Let us also note that the number of Newton iterations performed in
Algorithm \ref{alg_ip1} can be determined via Lemma \ref{l_anal}
(i). Moreover, if in Algorithm \ref{alg_ip2} we choose $\epsilon_V =
\delta_*/2$ and $\tau = 2 c /(2c +1)$ we need only one Newton
iteration at the inner stage. It follows that for this particular
choice for $\epsilon_V$ and $\tau$ the total number of Newton
iterations  of the algorithm   is given by the number of outer
iterations, i.e. the algorithm terminates in polynomial-time, within
${\cal O} \big ( \frac{1}{\log (\tau^{-1})}  \log (N_\phi
t^0/\epsilon) \big )$ iterations. This choice is made only for a
worst-case complexity analysis. In a practical implementation one
may choose larger values using heuristic considerations.


\section{Applications with Separable Structure}
\label{s_apl}

In this section we briefly  discuss some  of the applications to
which our method can be applied: distributed model predictive
control and network optimization. Note that  for these applications
our  Assumption \ref{ass1} holds.

\subsection{Distributed Model Predictive Control}

A first application that we will discuss  here is the control of
large-scale systems with interacting subsystem dynamics. A
distributed model predictive control (MPC) framework is appealing in
this context since this framework allows us to design local
subsystem-base controllers that take care of the interactions
between different subsystems and physical constraints. We assume
that the overall system model can be decomposed  into $N$
appropriate subsystem models:
\begin{align*}
x^i(k+1) = \sum_{j \in {\mathcal N}(i)} A_{ij}  x^j(k) + B_{ij}
u^j(k) \;\;  \forall i=1 \cdots N,
\end{align*}
 where ${\mathcal N}(i)$ denotes the set of
subsystems that interact with the $i$th subsystem, including itself.
The control and state sequence  must satisfy local constraints:
$x^i(k) \in \Omega_i$ and $u^i(k) \in U_i$  for all $i$ and  $k \geq
0$, where the sets $\Omega_i$ and $U_i$ are usually convex compact
sets with the origin in their interior (in general box constraints).
Performance is expressed via a stage cost, which we assume to have
the following form: $ \sum_{i=1}^N \ell_i(x^i,u^i)$, where usually
$\ell_i$ is a convex quadratic function, but not strictly convex in
$(x^i,u^i)$. Let $N_p$ denote the prediction horizon. In MPC we must
solve at each step $k$, given $x^i(k) = x^i$, an optimal control
problem of the following form \cite{MayRaw:00}:
\begin{align}
\label{dmpc}
 \!\!\min_{x_l^i,u_l^i} \big \{\!\! \sum_{l=0}^{N_p-1} \! \sum_{i=1}^N
\ell_i(x^i_l,u^i_l):  x_0^i\!=\!x^i, \ x_{l+1}^i \!=\!\sum_{j \in
{\mathcal N}(i)} A_{ij} x^j_l \!+\! B_{ij} u^j_l, \ x_l^i \in
\Omega_i, \ u_l^i \in U_i \ \forall  l,i \!\big \}.
\end{align}

 A similar formulation of distributed MPC for coupled
linear subsystems with decoupled costs was given in
\cite{VenRaw:07}, but without  state constraints. In
\cite{VenRaw:07}, the authors proposed to solve the optimization
problem \eqref{dmpc}  in a decentralized fashion, using the Jacobi
algorithm \cite{BerTsi:89}. But, there is no theoretical guarantee
of the Jacobi algorithm about how good  the approximation to the
optimum is after a number of iterations and moreover one needs
strictly convex functions $f_i$ to prove asymptotic convergence to
the optimum.

 Let us introduce $\xb_i
=[x_0^i \cdots x_N^i \ u_0^i \cdots u_{N-1}^i], X_i = \Omega_i^{N+1}
\times U_i^N$ and the self-concordant functions $f_i(\xb^i) =
\sum_{l=0}^{N_p-1} \ell_i(x^i_l,u^i_l)$  (recall  that $\ell_i$ are
assumed to be convex quadratic). The control problem \eqref{dmpc}
can be recast then as a separable convex program
\eqref{scp_f}--\eqref{scp_B}, where the matrices  $A_i$'s  and
$B_i$'s are defined appropriately, depending on the structure of the
matrices $A_{ij}$ and $B_{ij}$. In conclusion, Assumption \ref{ass1}
holds for this control problem so that our decomposition method can
be applied.


\subsection{Network Optimization}
\label{sec_no}
 Network optimization  furnishes another area in which
our algorithm leads to a new method of solution.  The  optimization
problem for routing in  telecommunication  data networks has the
following form \cite{XiaBoy:04,Lem:06}:
\begin{align}
\label{general_dcop}
 \min_{x_i \in [0, \ \bar x_i], y_j \in [0, \ d_j]}  \big \{
\sum_{j=1}^{n} f_j(y_j)+ \sum_{i=1}^{N} \langle c_i ,x_i \rangle: A
x_i =a_i, \; \sum_{i=1}^{N}  x_i =y \big \},
\end{align}
 where we  consider a multicommodity flow model with
$N$ commodities and $n$ links.  The matrix $A \in \rset^{m \times
n}$ is the node-link incidence matrix  representing the network
topology with entries $\{-1, 0, 1\}$. One of the most common cost
functions used in the communication network literature is the total
delay function \cite{XiaBoy:04,Lem:06}: $f_j(y_j) =
\frac{y_j}{d_j-y_j}$.
\begin{corollary}
Each function $f_j \in {\cal C}^3 \big ( [0, \ d_j) \big )$ is
convex  and $f_j$ is $3$-compatible with the self-concordant barrier
$\phi_j(y_j)=-\log(y_j(d_j-y_j))$ on the interval $(0, \ d_j)$.
\end{corollary}
\proof Note that the inequality \eqref{l_box} holds for all $y_j \in
(0, \ d_j)$ and $ h \in \rset$. Indeed,
\begin{align*}
|\nabla^3 f_j(y_j) | = 3 \nabla^2 f_j(y_j)  \sqrt{1/(d_j-y_j)^2}
\leq 3 \nabla^2 f_j(y_j)  \sqrt{1/(d_j-y_j)^2+ 1/y_j^2}. \qquad
\qed
\end{align*}

 Therefore, we can solve this  network optimization problem with
our method. Note that the standard dual function
 $d_0$ is not differentiable since it is the sum of a differentiable
  function (corresponding to the variable $y$) and a polyhedral
  function (corresponding to the variable $x$).  In \cite{Lem:06} a
  bundle-type algorithm is developed for maximizing the non-smooth function $d_0$, in
\cite{XiaBoy:04} the dual subgradient method is applied for
maximizing  $d_0$, while in \cite{BerTsi:89,KonLeo:96}  alternating
direction methods were proposed.

\subsection{Preliminary Numerical Results}
\label{simulations}

We illustrate the efficiency of our method in Table 1  on a random
set of problems of the form   \eqref{general_dcop}  and
\eqref{dmpc}, i.e. with  total delay (first half table) and
quadratic (second half) objective function, respectively.  For the
quadratic test problems we generate randomly the Hessian such that
it is positive semidefinite of the form $Q_i^T Q_i$, where $Q_i$ are
full row rank matrices. Here, the sets $X_i$ are assumed to have the
form $[0,\ u_i] \subseteq \rset^{n_1}$, i.e. $n_i= n_1$ and also
$m_i = m_1$ for all $i$. Note that for these type of problems the
barrier parameters $N_i= 2 n_1$  and $\alpha \leq 8$ and thus $c =
c_1 + c_2 \sqrt{n_1}$, for appropriate  $c_i > 0$ derived from Lemma
\ref{l_sscfd}. In our simulations we take $\tau = 0.85$, although  a
better tuning of this parameter will lead to less number of
iterations. Complexity  bounds for long step updates can also be
derived using similar arguments as those given in the present paper
for the short step method (see also Section 3.2.6 in
\cite{NesNem:94}). For all test problems the coupling constraints
have the form $\sum_i x_i = b$, so that the total number of
constraints is equal to $N m_1 + n_1$.

 In the table
we display the CPU time (seconds) and   the number of calls of the
dual function (i.e. the total number of outer and inner iterations)
for our dual interior-point algorithm (DIP) and an algorithm  based
on alternating direction method \cite{KonLeo:96} (ADI) for different
values of $m_1, n_1, N$ and fixed accuracy $\epsilon= 10^{-4}$. For
two problems the ADI algorithm did not produce the result after
running one day. All codes are implemented in Matlab version 7.1 on
a Linux operating system for both methods. The computational time
can be considerably reduced, e.g. by treating sparsity using more
efficient techniques as explained in Section \ref{s_PI} and
programming the algorithm
 in C. There are primal-dual interior-point methods that  treat
sparsity very efficiently but most of them specialized to
block-angular linear programs  \cite{GonSar:03}.  For different data
but with the same dimension and structure we observed that the
number of iterations does not vary much.

\begin{table}[t]
\label{table_sim}
\begin{center}
   \begin{tabular}{| c | c | c || c | c || c | c| }
     \multicolumn{3}{c}{} &   \multicolumn{2}{c}{DIP}         & \multicolumn{2}{c}{ADI}\\
     \hline
     $m_1$  &$n_1$  &  $N$    &   CPU   &   fct. eval.   &   CPU    &    fct. eval.      \\   \hline
       20   &  50   &  10     &   7.85  &    58          &    61.51  &     283          \\   \hline
       25   &  50   &  20     &   16.88 &    82          &  145.11   &     507           \\   \hline
       50   &  150  &   50    &  209.91 &    185         &  4621.42&     1451          \\    \hline
       80   &  250  &  100    & 1679.81 &    255         &  16548.23&     1748           \\   \hline
       170  &  500  &  100    & 10269.12&    367         &      *   &     *           \\  \hline \hline
       20   &   50  &   30    &   19.02 &     95         &   182.27  &     542           \\   \hline
       40   &  100  &  40     &   143.7  &    152         &   3043.67 &     1321          \\    \hline
       60   &  150  &   50    &   229.32&    217         &  10125.42&     2546          \\    \hline
       90   &  250  &  100    &  2046.09&    325         &  32940.67&     3816          \\   \hline
      100   &  300  &  120    &  4970.52&    418         &  *       &      *          \\   \hline
   \hline
   \end{tabular}
\end{center}
\caption{Computational results for network problems
\eqref{general_dcop} (first half) and quadratic problems
\eqref{dmpc} (second half) using DIP and ADI algorithms.}
\end{table}



\section{Conclusions}
A new decomposition method in convex programming is  developed in
 this paper using   dual decomposition   and interior-point framework.  Our method
 combines the fast local convergence rates of the Newton  method with the
 efficiency of structural optimization for solving
 separable convex programs. Although our  algorithm
 resembles  augmented Lagrangian methods, it  differs both in the
 computational steps and in the choice  of the parameters. Contrary to most augmented Lagrangian
 methods that use gradient based directions to update the Lagrange
 multipliers, our method  uses Newton directions and thus the
 convergence rate of the proposed method is faster. The reason for
 this  lies  in the fact that by adding self-concordant barrier
 terms  to the standard Lagrangian we proved that under appropriate
 conditions the corresponding family of augmented dual functions is
 also self-concordant. Another appealing theoretical   advantage  of our interior-point Lagrangian decomposition
method  is that it is  fully automatic, i.e. the parameters  of the
scheme are chosen  as in the
  path-following methods, which are crucial for justifying its global  convergence
 and polynomial-time  complexity.


%

\setstretch{0.8}


\begin{thebibliography}{10}

\bibitem{XiaBoy:04}
Xiao, L., Johansson, M., and Boyd, S.,
\newblock {\em Simultaneous routing and resource allocation via dual
decomposition},
\newblock {IEEE Transactions on Communications}, Vol. 52, No. 7, pp. 1136--1144, 2004.


\bibitem{GonSar:03}
Gondzio, J., and Sarkissian, R.,
\newblock {\em Parallel interior point solver for structured linear
programs},
\newblock {Mathematical Programming}, Vol. 96, pp. 561--584, 2003.


\bibitem{VenRaw:07}
Venkat, A., Hiskens, I., Rawlings, J., and Wright, S.,
\newblock {\em Distributed {MPC} strategies with application to power system
  automatic generation control},
\newblock {IEEE Transactions on Control Systems Technology}, to appear,
  2007.


\bibitem{Nec:08}
Necoara, I. and  Suykens, J. A. K.,
\newblock {\em Application of a smoothing technique to decomposition in convex
  optimization},
\newblock {IEEE Transactions on Automatic Control}, Vol. 53, No 11, pp. 2674--2679, 2008.


\bibitem{Zha:05}
Zhao, G.,
\newblock {\em A Lagrangian dual method with self-concordant barriers for
  multi-stage stochastic convex programming},
\newblock {Mathematical Programming}, Vol. 102, pp. 1--24, 2005.


\bibitem{BerTsi:89}
Bertsekas, D. P., and Tsitsiklis, J. N.,
\newblock {\em Parallel and distributed computation: Numerical
Methods},
\newblock Prentice-Hall, Englewood Cliffs, NJ, 1989.


\bibitem{Lem:06}
Lemarechal, C.,  Ouorou, A., and Petrou, G.,
\newblock {\em A bundle-type algorithm for routing in telecommunication data
  networks},
\newblock { Computational Optimization and Applications}, 2008.


\bibitem{KonLeo:96}
Kontogiorgis, S., De Leone, R., and Meyer, R.,
\newblock {\em Alternating direction splittings for block angular parallel
  optimization},
\newblock  {Journal of Optimization Theory and Applications}, Vol. 90, No. 1, pp. 1--29,
  1996.


\bibitem{CheTeb:94}
Chen, G. and Teboulle, M.,
\newblock {\em A proximal-based decomposition method for convex minimization
  problems},
\newblock {Mathematical Programming}, Vol 64, pp. 81--101, 1994.



\bibitem{KojMeg:93}
Kojima, M., Megiddo, N., Mizuno, S., and Shindoh, S.,
\newblock {\em Horizontal and vertical decomposition in interior point methods for
  linear programs},
\newblock  Dept. of Mathematical and Computing Sciences Technical report, Tokyo
  Institute of Technology, 1993.

\bibitem{KorPot:91}
Kortanek, K. O., Potra, F., and Ye, Y.,
\newblock {\em On Some Efficient Interior-Point Methods for Nonlinear Convex
Programming},
\newblock {Linear Algebra and Its Applications}, Vol. 152, pp. 169-189, 1991.


\bibitem{Tse:92}
Tseng, P.,
\newblock {\em Global Linear Convergence of a Path-Following Algorithm for
 Some Monotone Variational Inequality Problems},
\newblock {Journal of Optimization Theory and Applications}, Vol. 75, No. 2, pp.
265-279, 1992.


\bibitem{Zhu:92}
Zhu, J.,
\newblock {\em A path following algorithm for a class of
convex programming problems},
\newblock  {Mathematical Methods of Operations Research}, Vol. 36,
 No. 4, pp. 359-377, 1992.


\bibitem{HegOsb:01}
Hegland, M., Osborne, M. R., and  Sun, J.,
\newblock {\em  Parallel Interior Point Schemes for Solving Multistage Convex
Programming},
\newblock {Annals of Operations Research},  Vol. 108, No 1--4, pp. 75-85, 2001.




\bibitem{MieMos:72}
Miele, A., Moseley, P. E., Levy, A. V., and Coggins, G. M,
\newblock {\em On the method of multipliers for mathematical programming problems},
\newblock  {Journal of Optimization Theory and Applications}, Vol. 10, pp. 1--33,
  1972.


\bibitem{NesNem:94}
Nesterov, Y., and Nemirovskii, A.,
\newblock {\em Interior Point Polynomial Algorithms in Convex
Programming},
\newblock Society for Industrial and Applied Mathematics (SIAM Studies in
  Applied Mathematics), Philadelphia, 1994.


\bibitem{Ren:01}
Renegar, J.,
\newblock{\em A Mathematical View of Interior-Point Methods for Convex
Optimization},
\newblock MPS-SIAM Series on Optimization, Philadelphia, 2001.



\bibitem{MayRaw:00}
Mayne D. Q., Rawlings J. B., Rao C. V., and  Scokaert P. O. M.,
\newblock {\em Constrained Model Predictive Control: Stability and Optimality},
\newblock  {Automatica},  Vol. 36,  No. 7,  pp.  789--814, 2000.

\end{thebibliography}
\end{document}